\documentclass[12pt,reqno]{amsart}
\pdfoutput=1
\usepackage[headings]{fullpage}
\usepackage{amssymb,bbm}
\usepackage{epstopdf}
\usepackage{graphicx}
\usepackage[utf8]{inputenc}
\usepackage[russian,english]{babel}
\usepackage{texdraw}
\usepackage{url}
\usepackage[all]{xy}
\usepackage{pb-diagram}
\usepackage{float}   
\usepackage[bookmarks=true,%
    colorlinks=true,%
    linkcolor=blue,%
    citecolor=blue,%
    filecolor=blue,%
    menucolor=blue,%
    urlcolor=blue,%
    breaklinks=true]{hyperref}

\theoremstyle{definition}

\def\BZ{\mathbbm Z}
\def\BQ{\mathbbm Q}
\def\BR{\mathbbm R}
\def\BC{\mathbbm C}
\def\BP{\mathbbm P}
\def\BH{\mathbbm H}

\def\B{\mathcal B}  \def\hB{\widehat\B}

\def\calL{\mathcal L}
\def\calO{\mathcal O}
\def\calT{\mathcal T}
\def\calH{\mathcal H}

\def\b{\beta}
 
\def\G{\Gamma}
\def\la{\langle}
\def\ra{\rangle}

\def\wh{\widehat}
\def\SL{\mathrm{SL}}
\def\PSL{\mathrm{PSL}}
\def\GL{\mathrm{GL}}
\def\PGL{\mathrm{PGL}}
\def\Sp{\mathrm{Sp}}

\def\D{\Delta}

\def\Vol{\mathrm{Vol}} 
\def\V{{\Vol_{\C}} }
\def\WW{\mathsf\L^2} 
\def\a{\alpha} 
\def\k{\kappa}
\def\z{\zeta}

\def\={\;=\;}
\def\+{\,+\,}
\def\m{\,-\,} 
\def\-{\,-\,}

\def\be{\begin{equation}}
\def\ee{\end{equation}}
\def\bea{\begin{align}}
\def\eea{\end{align}}

\def\sm{\smallsetminus}

\def\Li{\mathrm{Li}}
\newcommand{\mb}{\mathbf}

\def\m{\,-\,}
\def\R{\mathbb R}
\def\C{\mathbb C}
\def\Z{\mathbb Z}

\def\Q{\mathbb Q}

\def\th{\theta}
\def\O{\text O} 

\def\J{\mathbf J}
\def\bPhi{\mathbf\Phi}
\def\bp{\mb p}  \def\bq{\mb q}  \def\bz{\mb z} 
\def\sma#1#2#3#4{\bigl(\smallmatrix#1&#2\\#3&#4\endsmallmatrix\bigr)}

\def\Im{\mathrm{Im}}

\def\Qbar{\overline{\BQ}}
\def\g{\gamma}
\def\ssm{\sm}

\def\l{\lambda}

\def\ve{\varepsilon}

\def\sma#1#2#3#4{\bigl(\!\smallmatrix#1&#2\\#3&#4\endsmallmatrix\!\bigr)}  
\def\vsma#1#2#3#4{(\smallmatrix#1&#2\\#3&#4\endsmallmatrix)} 

\def\Lob{\text{\foreignlanguage{russian}{Л}}}

\def\fu{\mathfrak u}   
\def\fv{\mathfrak v}   
\def\fU{\mathfrak U}   
\def\L{\Lambda}

\makeatletter

\renewcommand\thepart{\@Roman\c@part}%
\renewcommand\part{%
   \if@noskipsec \leavevmode \fi
   \par
   \addvspace{6.7ex}%
   \@afterindentfalse
   \secdef\@part\@spart}
\def\@part[#1]#2{%
    \ifnum \c@secnumdepth >\m@ne
      \refstepcounter{part}%
      \addcontentsline{toc}{part}{Part~\thepart.\ #1}%
    \else
      \addcontentsline{toc}{part}{#1}%
    \fi
    {\parindent \z@ \raggedright
     \interlinepenalty \@M
     \normalfont
     \ifnum \c@secnumdepth >\m@ne
       \centering\large\scshape \partname~\thepart.%
       \hspace{1ex}%
     \fi%
     \large\scshape #2%
     \markboth{}{}\par}%
    \nobreak
    \vskip 4.7ex
    \@afterheading}
  \def\@spart#1{
  \refstepcounter{part}%
  \addcontentsline{toc}{part}{#1}%
    {\parindent \z@ \raggedright
     \interlinepenalty \@M
     \normalfont
     \centering\large\scshape #1\par}%
     \nobreak
     \vskip 4.7ex
     \@afterheading}
\renewcommand*\l@part[2]{%
  \ifnum \c@tocdepth >-2\relax
    \addpenalty\@secpenalty
    \addvspace{0.75em \@plus\p@}%
    \begingroup
      \parindent \z@ \rightskip \@pnumwidth
      \parfillskip -\@pnumwidth
      {\leavevmode
       \normalsize \bfseries #1\hfil \hb@xt@\@pnumwidth{\hss #2}}\par
       \nobreak
       \if@compatibility
         \global\@nobreaktrue
         \everypar{\global\@nobreakfalse\everypar{}}%
      \fi
    \endgroup
  \fi}

\def\l@subsection{\@tocline{2}{0pt}{2pc}{6pc}{}}
\makeatother

\begin{document}
\title[Hyperbolic 3-manifolds, the Bloch group, and the work of Walter
Neumann]{Hyperbolic 3-manifolds, the Bloch group, \protect
  \\ and the work of Walter Neumann}
\author{Stavros Garoufalidis}
\address{
  International Center for Mathematics, Department of Mathematics \\
  Southern University of Science and Technology \\
  Shenzhen, China \newline
  {\tt \url{http://people.mpim-bonn.mpg.de/stavros}}}
\email{stavros@mpim-bonn.mpg.de}
\author{Don Zagier}
\address{Max Planck Institute for Mathematics \\
  Bonn, Germany and
International Centre for Theoretical Physics \\ Trieste, Italy}
\email{dbz@mpim-bonn.mpg.de}


% ...

\maketitle

{\footnotesize
\tableofcontents
}


Walter Neumann---henceforth simply Walter---has made huge contributions to
topology, singularity theory, and number theory.  His work on the combinatorics
of triangulations of 3-manifolds and its interactions with algebraic $K$-theory
certainly can be counted as among the most far-reaching and beautiful of these.
In this short note we will try to recount some of the high points of this work,
taking as our starting point the papers~\cite{NZ} (by Walter and the second
author) and~\cite{N:combinatorics} (by him alone). In the first sections we
describe the contents of these two papers and of related later work of 
Walter, both alone and with other coauthors, in some detail, 
while the final sections will present further development of some of these 
themes in later work of other people, including ourselves.


\section{Ideal triangulations and the gluing equations}
\label{sec.2}

The starting point of the paper~\cite{NZ} was Thurston's amazing insight in the 1980's
that all 
3-dimensional manifolds should be canonically divisible into pieces having a well-defined 
geometric structure of one of 8 types, the most important of which is the hyperbolic one. 
In conjunction with the famous Mostow ridigity theorem, this means that 3-dimensional
topology becomes a part, first of differential geometry, and then of algebraic number
theory, something that is not at all the case in other dimensions.  The main class is
that of oriented hyperbolic 3-manifolds, which have a riemannian metric with constant
negative curvature that can be normalized to~$-1$ and hence are locally isometric to
hyperbolic 3-space~$\BH^3$.  Of particular interest is the
{\it volume spectrum}, the set of volumes of all complete hyperbolic 3-manifolds of finite
volume (in which case they are either compact or the union of a compact part and a finite
number of ``cusps'' diffeomorphic to the product of a half-line and a torus).  These
volumes have both striking number-theoretical properties (they belong to the image under
the regulator map of the Bloch group, as discussed in Section~\ref{sec.4}) and striking
metric properties (they form a countable well-ordered subset of~$\R_{>0}$, as discussed in
Section~\ref{sec.3}), and the primary goal of the paper~\cite{NZ} was to understand them
as thoroughly as possible.

In this section we discuss ideal triangulations and their NZ-equations in some detail. 
Ideal triangulations of 3-manifolds with torus boundary components were introduced by
Thurston~\cite{Thurston} as a convenient way to describe and effectively
compute~\cite{snappy} complete hyperbolic structures on 3-manifolds.
Recall that in hyperbolic geometry an {\it ideal tetrahedron} is the convex
hull of four points in the boundary $\partial (\BH^3)\cong\BP^1(\BC)$ of 3-dimensional
hyperbolic space $\BH^3$.  The (orientation preserving) isometry group of~$\BH^3$ is the
group~$\PSL_2(\BC)$, 
acting on the boundary by fractional linear transformations of~$\BP^1(\BC)$, and since under 
this group any four distinct points can be put in standard position~$(0,1,\infty,z)$ for some 
$z\in\BC\ssm\{0,1\}$ (the cross-ratio of the four points), any (oriented) ideal tetrahedron is 
the convex hull $\D(z)$ for some complex number~$z~\in \BC\ssm\{0,1\}$, called the
{\it shape parameter} 
of the tetrahedron. This number is not quite unique because of the choice of which three
vertices
we send to 0, 1 and~$\infty$, meaning that the oriented tetrahedra $\D(z)$, $\D(z')$ and
$\D(z'')$,
where $z'=1/(1-z)$ and $z''=1-1/z$, are isometric. In this way one attaches a shape parameter 
$z$, $z'$ or $z''$ to each pair of opposite edges of a given oriented ideal tetrahedron.
Note that the shape of a tetrahedron is an arbitrary complex number not equal to 0 or 1,
and whether it has positive, zero or negative imaginary part is irrelevant to our discussion.

When $\calT$ is an ideal triangulation of a hyperbolic 3-manifold $M$ with $N$ tetrahedra
with shape parameters $z_1,\dots,z_N$, then we get one polynomial equation
(``gluing equation'') 
for each edge.  Specifically, the shape parameters at that edge of all tetrahedra that are
incident with it must clearly have arguments that add up to~$2\pi$ (because otherwise the 
metric would not be smooth along that edge), but in fact the shape parameters themselves
have product~$+1$ by an easy argument.  Since each of the possible shape parameters 
$z,\,z',\,z''$ of~$\D(z)$ belongs to the multiplicative
group~$\big\langle z,1-z,-1\big\rangle$, this equation for each edge~$e_i$  has the form 
\be
\label{edge}
\pm\prod_{j=1}^N z_j^{R'_{ij}} (1-z_j)^{R''_{ij}} = 1\,.
\ee

Since it is easily seen that the number of edges in the triangulation is the same as 
the number~$N$ of simplices, this gives us $N$ polynomial equations among the $N$ complex
numbers $z_1,\dots,z_N$.  The obvious thought is that this leads to a 0-dimensional 
moduli space and explains the rigidity, but this is wrong since rigidity only applies 
to complete hyperbolic structures, and in fact the $N$ edge relations are never
multiplicatively
independent. An obvious example is that their product is always~1, and this is the
only dependence if the boundary of the 3-manifold is a torus (often called a cusp),
but in general there are $h$ multiplicatively independent relations among the $N$
equations~\eqref{edge}, where $h$ is the number of cusps of the 3-manifold (assuming that
all boundary components are tori). Thus the true expected 
dimension of the moduli space of hyperbolic structures is in fact~$h$.  That the 
dimension really is~$h$ is an important theorem of Thurston (\S5 of~\cite{Thurston}) and 
was given a new and simpler proof in~\cite{NZ} using the algebraic structure of the 
gluing equations.
The 0-dimensional moduli space (rigidity) arises when we require the hyperbolic structure
given by the shape parameters $z_i$ to be complete, because this entails two further
independent relations at each cusp.  Specifically, each ``peripheral curve'' (meaning an 
isotopy class of curves on the torus cross-section of one of the cusps) gives a relation,  
so if one chooses a meridian and a longitude $(\mu_i,\lambda_i)$ at each cusp 
$i=1,\dots,h$, one obtains $2h$ further equations
\be \label{ML}
\pm\prod_{j=1}^N z_j^{M'_{ij}} (1-z_j)^{M''_{ij}} \= 1, \qquad
\pm\prod_{j=1}^N z_j^{L'_{ij}} (1-z_j)^{L''_{ij}} \= 1, \qquad (i=1,\dots,h) \,.  
\ee
Thus the full set of gluing equations is described by an $(N+2h) \times 2N$ matrix 
$U=\Bigl(\!\begin{smallmatrix} R'&R''\\M'&M''\\L'&L''\end{smallmatrix}\!\Bigr)$. 

It was shown in~\cite{NZ} that the matrix~$U$ has some key symplectic properties, of which
the most important (Theorem~2.2) says that
\be\label{Th2.2} 
U J_{2N}U^t \= 
\begin{pmatrix} 0_{N,N}& 0_{2h,N} \\ \phantom2 0_{N,2h}  & 2 J_{2h} \end{pmatrix}\,,
\ee
where $J_{2n}:=\bigl(\!\begin{smallmatrix} 0&-1_n\\1_n&0\end{smallmatrix}\!\bigr)$, while
the others say that the $N\times2N$ matrix $R=\bigl(R'\,R''\bigr)$ has rank $N-h$,
the full matrix~$U$
has rank $N+h$, and that the space $[U]\subset\R^{2N}$ spanned by the rows of~$U$ is the 
orthogonal complement of~$[R]$ with respect to the symplectic structure~$J_{2N}$ (Prop.~2.3).
The rank statement was then used to prove that the dimension of the above-mentioned
deformation space of (non-complete) hyperbolic structures is at least~$h$, and then a
further argument using Mostow rigidity showed that it is exactly~$h$.  This gives rise in
the 1-cusp case to a polynomial $P(m,\ell)=0$, where $m^2$ and $\ell^2$ denote the
left-hand sides of the two equations~\eqref{ML} and $P$ is a certain polynomial (a factor
of what is now called the $A$-polynomial) that was calculated explicitly in~\cite{NZ} for 
the simplest hyperbolic knot~$4_1$ (figure~8).  In the general case the deformation space
is a component of the ``character variety'' as introduced and studied in~\cite{CCGLS} and
studied later by Walter and his students Abhijit Champanerkar~\cite{Champanerkar} and
Stefan Tillmann~\cite{Tillmann:bs}.

To get a complete manifold with $h$ cusps, we impose $h$ further gluing conditions, namely the
equations~\eqref{edge} and the first of each of the equations~\eqref{ML}. (But in fact these
$h$ meridian equations together with the edge equations imply the longitude equations, so in
the end all relations~\eqref{edge} and~\eqref{ML} hold.)  Now rigidity applies and there are 
no deformations.  But there is another process, with a lot more freedom, to obtain rigid
complete hyperbolic structures from the original cusped manifold, namely to do a Dehn 
surgery at some or all of the cusps.  Specifically, when we do a $(p,q)$-surgery at a
cusp (meaning that we truncate the \hbox{3-manifold} at the torus boundary and glue on a solid
torus in such a way as to kill the homotopy class of $p$ times a chosen meridian times $q$
times a chosen longitude), then we impose the gluing equation $m^{2p}\ell^{2q}=1$, where
$m^2$ and $\ell^2$ as above are the left-hand sides of one of the pairs of
equations~\eqref{ML}.  In the next section, still following~\cite{NZ}, we discuss how the
volumes behave under this process.


\section{Volumes and Dehn surgeries}
\label{sec.3}
 
As already stated, the original main purpose of~\cite{NZ} was to study the arithmetic of 
the set of volumes of all hyperbolic 3-manifolds. It had been observed by Thurston, using
earlier work of J\o rgensen, that this volume spectrum is a well-ordered subset of the
positive reals.
In other words, there is a smallest volume (which is known), a second smallest, a third 
smallest, \dots, then a smallest limit point, a second smallest limit point, \dots, then
limits of 
these, and so on.  The proof of well-orderedness shows that these simple and higher-order
limit
points arise by ``closing up''  one or more of the cusps of a non-compact hyperbolic
3-manifold~$M$ by Dehn surgeries to obtain
a countable collection of manifolds with fewer cusps whose volumes tend from below to 
that of~$M$.  The object of~\cite{NZ} was to study the speed with which these volumes converge.

Before studying the effect of surgeries, we must understand the volume of a single hyperbolic 
3-manifold~$M$.  Clearly it can be given as the sum of the volumes of the tetrahedra of any
ideal triangulation, so the
first step is to understand these.  An (oriented and non-degenerate) ideal tetrahedron can be 
parametrized either by a shape parameter~$z$ in the complex upper half-plane, as explained in
the preceding section, or, in case~$z$ is in the upper half plane,
by the three angles $\a,\,\b,\,\g$ (positive and
with~$\a+\b+\g=\pi$) of the Euclidean triangle that one ``sees'' by looking at the
tetrahedron from any of its four cusps. When $\Im(z)>0$, the angles and the shapes are
related by $(\a,\b,\g)=(\arg(z),\arg(z'),\arg(z''))$ with $z'$ and~$z''$ as above.

According to Chapter~7 of Thurston's notes, written by Milnor, the volume of this
tetrahedron when $\Im(z)>0$ is given in terms of these two parametrizations by the two formulas
\be
\label{Kummer}
\Vol(\D(z)) \=  D(z) \;=\; \Lob(\a)+\Lob(\b)+\Lob(\g)\,,
\ee
the equality of the two being an identity due to Kummer.  (Both formulas are actually true
for all $z\ne0,\,1$, but the sum of $\a$, $\b$ and~$\g$ is $-\pi$ if $z$ is in the lower half-plane.)
Here $D(z)$ and $\Lob(\th)$ are the Bloch-Wigner dilogarithm (cf.~Section~\ref{sec.4}) and the 
Lobachevsky function, defined respectively by 
\be
\label{defD}
\begin{aligned}  D(z) &\=\Im\bigl(\Li_2(z)+\log|z|\log(1-z)\bigr)\,, \\
  \Lob(\th) &\= -\int_0^{\th} \log|2\,\sin t|\,dt
  \=\frac12\,\sum_{n=1}^\infty\frac{\sin(2n\th)}{n^2}
 \= \frac12\,D\bigl(e^{2i\th}\bigr)\;.
\end{aligned}
\ee
Thus the volume of a hyperbolic manifold $M$ triangulated by $N$ ideal tetrahedra with shape 
parameters~$z_i$ is given in terms of the dilogarithm function by
\be
\label{volM}  \Vol(M) \= \sum_{j=1}^N D(z_j)\;.
\ee

Now let $M$ be a 3-manifold with $h$ cusps.  It has a unique complete hyperbolic structure as
a quotient $\BH^3/\G$ for some lattice $\G\subset\PSL_2(\C)$ isomorphic to~$\pi_1(M)$, the
corresponding shape parameters $z^0=(z^0_1,\dots,z^0_N)$ being a solution of all $N+2h$
equations~\eqref{edge} and~\eqref{ML}. A small deformation of this hyperbolic structure
will have a nearby shape parameter vector $z=(z_1,\dots,z_N)$ satisfying only~\eqref{edge}.
By the result of Thurston mentioned earlier, the space of all such deformations has the
structure of a smooth complex manifold of dimension~$h$, so is isomorphic to a small
neighborhood $\fU$ of~0 in~$\C^h$. Thus each $z_i$ depends holomorphically on $\fu\in\fU$
and~$z_i(0)=z_i^0$.  To make this more visible, and to choose nice coordinates $\fu_i$
on~$\fU$, we have to look more closely at the structure of the cusps on~$M$ and at Dehn
surgeries.

Each cusp has the structure $[0,\infty)\times T^2$, where $T^2$ is a totally geodesically
embedded torus in~$M$ and has a flat Euclidean metric, unique up to homothety, induced by
the hyperbolic metric on~$M$.
More concretely, the cusps of~$M=\BH^3/\G$ are indexed by $P\in\BP^1(\C)$ whose
stabilizer~$\G_P$ is free
abelian of rank~2.  After conjugation we can place $P$ at~$\infty$, in which case 
$\G_\infty$ has the form
$\pm\Bigl(\begin{matrix} 1&\L\\0&1\end{matrix}\Bigr)$ for some lattice~$\L\subset\C$, 
unique up to homothety, and then we can identify~$\L$ with $H_1(T^2)$ and $T^2$ with~$\C/\L$. 
There is a canonical quadratic form on~$\L$ defined as the square of the length of a 
vector divided by the volume of~$T^2$.  If we choose an oriented basis $(\mu,\l)$ of 
$H_1(T^2)$ (``meridian'' and ``longitude'') to identify $\L$ with~$\Z^2$, then this quadratic 
form is given by $Q(p,q)=|p\tau+q|^2/\Im(\tau)$ if we rescale~$\L$ after homothety to 
be~$\Z\tau+\Z$ for some~$\tau$ in the complex upper half-plane.  

On the other hand, each element of $\L=H_1(T^2)$ can be identified with an isotopy class of
closed curves on~$T^2$. Doing a $(p,q)$ Dehn surgery at a cusp, once the basis of~$H_1(T^2)$ 
has been chosen, means removing $T^2\times[0,\infty)$ from~$M$ and replacing it by a solid 
torus in such way that the curve on~$T^2$ corresponding to the class $p\mu+q\l$ bounds 
in the solid torus. Here we assume that the integers $p$ and~$q$ are coprime, but we also 
allow~``$\infty$'' as a value for~$(p,q)$, meaning that we leave this cusp untouched.  The 
Dehn surgeries on~$M$ are then described by tuples 
$\k=((p_1,q_1),\dots,(p_h,q_h))\in(\Z^2\cup\{\infty\})^h$.  Thurston's theorem~\cite{Thurston}
tells us that if all $(p_i,q_i)$ are near enough to~$\infty$ (meaning that $(p_i,q_i)=\infty$ or 
$p_i^2+q_i^2$ is large), the surgered manifold~$M_\k$ is hyperbolic, in which case 
both the shape parameters $z_1,\dots,z_N$ and the volume~$\Vol(M_\k)$ become functions 
of~$\k$.  Explicitly, the deformed shape parameters 
$\bz(\k)=\bz(\bp,\bq)$ that tend to~$\bz^0$ as all the pairs $(p_i,q_i)$ tend to infinity 
are given by adding to the original gluing equations~\eqref{edge} the $h$ new equations given
by the product of the $p_i$-th power of the first expression by the $q_i$-th power of the 
second one in~\eqref{ML}, and then $\Vol(M_\k)$ is given by~\eqref{volM} with $z_j$
replaced by~$z_j(\k)$.

We now get two real numbers at each cusp: the value $Q_i(p_i,q_i)$ of the quadratic form
corresponding to that cusp at the pair $(p_i,q_i)$ (with the
convention~$Q_i(\infty)=\infty$) and the length
$L_i=L_i(\k)$ of the short geodesic on~$M_\k$ which is the core of the solid torus added
by the Dehn surgery
(or~0 if the $i$th cusp has not been surgered). They are related by
\be
\label{systole}
L_i \= \frac{2\pi}{Q_i(p_i,q_i)} \+ \O\biggl(\sum_{i=1}^h\frac1{p_i^4+q_i^4}\biggr)
\ee
(\cite{NZ}, Prop.~4.3) as all~$\k_i=(p_i,q_i)$ tend to infinity.  The main volume result
of~\cite{NZ}, proved by using~\eqref{volM} and analyzing the changes of the dilogarithms
under small changes of the~$z$'s, is then given by the pair of asymptotic formulas 
\begin{subequations}
  \label{XXX}
  \begin{align}\label{VolDehn1}
\Vol(M_\k) 
&\= \Vol(M) \;-\;\sum_{i=1}^h\biggl(\frac{\pi^2}{Q_i(p_i,q_i)}
\+ \O\Bigl(\frac1{p_i^4+q_i^4}\Bigr)\biggr) \\
\label{VolDehn2}  
&\= \Vol(M) \;-\;\sum_{i=1}^h\Bigl(\frac{\pi L_i}2 \+ \O\bigl(L_i^2\bigr)\Bigr)
\end{align} \end{subequations}
(Theorems 1A and~1B in~\cite{NZ}), which are equivalent to one another by virtue
of~\eqref{systole}. These volumes, as $\k=(\bp,\bq)$ ranges over all $h$-tuples of
sufficiently large pairs of coprime integers or the symbol~$\infty$ meaning unsurgered,
all belong to the hyperbolic volume spectrum, and equation~\eqref{VolDehn1} has as an
immediate corollary a description of the local structure of this volume spectrum near its
limit point, because the asymptotics of the number of lattice points, or of primitive
lattice points, in a large ellipse, is well-known.  The precise asymptotic statement, 
which we will not repeat here, is formulated explicitly as a Corollary to Theorem~1A
in~\cite{NZ}.

All of this is only for integral (and coprime) values of~$p_i$ and~$q_i$.  However, as is 
explained in~\cite{NZ} in detail, the shape parameters $z_j(\k)$ and the lengths $L_i(\k)$ 
are defined for $p$ and~$q$ real rather than just integral and coprime in pairs, and
equations~\eqref{systole} and~\eqref{XXX} still remains true.  The only point is that in
the definition of the $z_j$'s we had to take the $p_i$th and $q_i$th powers of the
equations in~\eqref{ML}, and one cannot in general take real powers of complex numbers
in a well-defined way, but since the left-hand sides of the expressions in~\eqref{ML}
are near to~1 for small deformations of the original value $\bz=\bz^0$ and since a
complex number near~1 has a well-defined logarithm near~0, there is no problem.  When the
$p_i$ and~$q_i$ are not integral, we are no longer ``filling in'' the cusp by gluing on
a solid torus, but are simply changing the hyperbolic structure on the original
open topological manifold~$M$, with the new  hyperbolic structures in general being
incomplete. 
If we now define $2h$ complex numbers $\fu=(\fu_1,\dots,\fu_h)$ and
$\fv=(\fv_1,\dots,\fv_h)$ by
$$ 
\fu_i\,=\,\sum_{j=1}^N\biggl(M'_{ij}\,\log\frac{z_j}{z_j^0}
\+M''_{ij}\,\log\frac{1-z_j}{1-z_j^0}\biggr),\quad
\fv_i\,=\,\sum_{j=1}^N\biggl(L'_{ij}\,
\log\frac{z_j}{z_j^0}\+L''_{ij}\,\log\frac{1-z_j}{1-z_j^0}\biggr),
$$
then from the symplectic properties of the gluing equations it follows
that~\cite[Lem.4.1]{N:combinatorics} 
$$ 
p_i\fu_i \+ q_i\fv_i\= 2\pi i\qquad(i=1,\dots,h)\,, 
$$
and we can take $\fu$ as canonical coordinates for the above-mentioned neighborhood $\fU$ 
of $0\in\C^h$, 
in which case each $\fv_i$ becomes an odd power series in the $\fu$'s with linear term 
$\tau_i\fu_i$
and we can write $L_i(\fu)$ and $\Vol(\fu)$ instead of $L_i(\k)$ and~$\Vol(M_\k)$.  
We should mention
that $\fu$ and~$\fv$ can be defined invariantly, without using any triangulation, 
as follows: the deformed
hyperbolic structure on~$M$ corresponds to a homomorphism $\rho:\G\to\PSL_2(\C)$ near to 
the inclusion
map, and then $\fu$ and~$\fv$ are simply the logarithms of the ratios of the eigenvalues of 
the images under~$\rho$ of the meridians and longitudes, respectively.

In terms of the new coordinates, equation~\eqref{VolDehn2} becomes
\be
\label{VolDehn3} 
\Vol(\fu) \= \Vol(M) \;-\;\frac\pi 2\,\sum_{i=1}^h L_i(\fu) \+ \ve(\fu),  
\ee
with $\ve(\fu)=\O(||\fu||^4)$ as $\fu$ tends to~0 in~$\C^h$. 
Theorem~2 of~\cite{NZ} was the statement that the function $\ve(\fu)$ defined 
by~\eqref{VolDehn3}
is harmonic, and hence is the real part of a holomorphic function $f(\fu)$ near~0 (uniquely
determined if we fix~$f(0)=0$).  
Theeorem~3, proved using equation~\eqref{Th2.2}, said that $\partial\fv_i/\partial\fu_j$ 
is symmetric in~$i$ and $j$, which implies that there is a single function $\Phi(\fu)$ with
$\partial\Phi/\partial\fu_i=2\fv_i$ for all~$i$, and that $f$ is given in terms of $\Phi$ by 
$4f=\Phi-\fu\cdot\fv$ (or equivalently by $-8f=(E-2)\Phi$, where
$E=\sum\fu_i\,\partial/\partial\fu_i$
is the Euler operator), so that the volume correction~$\ve(\fu)$ in~\eqref{VolDehn3} is
given by
\be
\label{VolDehn4}
\ve(\fu) \= \Im\bigl(f(\fu)\bigr),
\qquad f(\fu)\,:=\,\frac14\,\int_0^\fu
\biggl(\,\sum_{i=1}^h(\fv_i\,d\fu_i\,-\,\fu_i\,d\fv_i)\biggr)\,.
\ee
The function $f$ is now often called the Neumann-Zagier potential function, although
this name was used in the original paper for~$\Phi$ instead.  It should perhaps also be
mentioned that simpler proofs of the last results described could probably have been
obtained by using the second rather than the first volume formula in~\eqref{Kummer}.

There is one more important point about volumes. Another insight by Thurston was that the
volume of a hyperbolic 3-manifold, which is a positive real number, is  actually in a
natural way the imaginary part of a {\it complexified volume} whose real part is the
Chern-Simons invariant, an important topological 
invariant taking values in~$\R/4\pi^2\Z$ whose definition we omit here.
It was  conjectured in~\cite{NZ}, and proved soon afterwards by Yoshida~\cite{Yosh},
that the above formulas remain true with the volumes replaced by their complexified 
versions, the functions~$L_i(\fu)$ also lifted suitably from~$\R$ to~$\C$, and $\ve(\fu)$
replaced by~$f(\fu)$.  Later, in~\cite{N:combinatorics}, Walter showed how to
lift~\eqref{volM} to an explicit and computable expression for the complexified volume
of~$M_{\bp,\bq}$ in terms of the complex dilogarithm. 


\section{The Bloch group and the extended Bloch group} 
\label{sec.4}

The Bloch group of a field is an analogue of its multiplicative group, but with the 
relation $[xy]=[x]+[y]$ satisfied by the logarithm function replaced by the functional 
equation of the dilogarithm.  In this section we recall its definition and the definition 
of the ``extended Bloch group'' that was introduced by Walter~\cite{N:extended} and further 
developed by Zickert and Goette~\cite{Goette,Zickert:extended}, and explain their connections 
with the volume and complexified volume.  The next section tells how these things relate
to the symplectic structure.  We should mention that parts of both sections have been 
transferred here from the arXiv version of~\cite{GZ:kashaev} and also edited somewhat for 
the purpose of the present exposition.

The dilogarithm function~$\Li_2(z)$, defined for~$|z|<1$ as $\sum_{n=1}^\infty z^n/n^2$
and then extended analytically to either the cut plane $\C\ssm[1,\infty)$ or to the universal
cover of~$\BP^1(\C)\ssm\{0,1,\infty\}$, satisfies a famous functional equation called the
5-term relation. This functional equation was discovered repeatedly during the 19th century 
and can be written in many equivalent forms, each saying that a sum of five dilogarithm 
values is a linear combination of products of simple logarithms.  
The function~$\Li_2$ is many-valued, but the modified dilogarithm~\eqref{defD} is a 
single-valued real analytic function from $\BP^1(\C)\ssm\{0,1,\infty\}$ to~$\R$ that 
extends continuously to all of~$\BP^1(\C)$ and satisfies ``clean'' versions of the 5-term 
relations with no logarithmic correction terms.  Since $D(z)$ also satisfies the two 
functional equations $D(1-z)=-D(z)=D(1/z)$ (implying that $D(z)=D(z')=D(z'')$ for the 
three shape parameters of an oriented ideal hyperbolic tetrahedron), this ``clean" 
functional equation still can be written in many different forms, one standard one being
$$ 
D(x) \+ D(y) \+ D\Bigl(\frac{1-x}{1-xy}\Bigr) \+ D(xy) \+ D\Bigl(\frac{1-y}{1-xy}\Bigr) \= 0  
$$
for $(x,y)\ne(1,1)$ in~$\C^2$.  Another nice version is the cyclic one 
$\sum_{i\!\pmod5} D(z_i)=0$ if $\{z_i\}_{i\in\Z}$ is a sequence of complex numbers 
satisfying $1-z_i=z_{i-1}z_{i+1}$ for all~$i$ (which implies by a short calculation that 
they have period~5). Yet another, with a clear interpretation in terms of 3-dimensional 
hyperbolic geometry, says that the signed sum of $D(r_i)$ is~0 if $r_1,\dots,r_5$
are the cross-ratios of the 5 subsets of cardinality~4 of a set of 5 distinct points
in~$\BP^1(\C)$. 

The five arguments $z_i$ of any version of the five-term relation satisfy 
$\sum(z_i)\wedge(1-z_i)=0$, where the sum is taken in the second exterior power of the
multiplicative
group of~$\C$. (For instance, for the ``cyclic version" above we have $\sum_i(z_i)\wedge(1-z_i)
= \sum_i(z_i)\wedge(z_{i-1}z_{i+1})=\sum_i\bigl((z_i)\wedge(z_{i-1})+(z_{i-1})\wedge(z_i))=0$.)
The Bloch group $\B(F)$ of an arbitrary field~$F$, introduced by Bloch~\cite{Bloch} in 1978, 
is motivated by this observation and is defined
as the quotient of the kernel of the map $d:\Z[F]\to\WW(F^\times)$ sending $[x]$ to 
$x\wedge(1-x)$ for $x\ne0,1$ (and to~0 for $x=0,1$) by the subgroup generated by the 5-term
relation of the dilogarithm. The precise definition varies slightly in the literature
because of 
delicate 2- and 3-torsion issues arising from the particular definition of the exterior square 
(for instance, does one require $x\wedge x=0$ for all $x$ or just
$x\wedge y=-y\wedge x\,$?), whether
one requires $d([0])$ and $d([1])$ to vanish or merely to be torsion, and the particular
version of 
the 5-term relation used.  We will gloss over this point for now, but will come back to it in
connection with the extended Bloch group.

From our point of view, the clearest motivation for the definition of the Bloch group is
the fact that the shape parameters~$\{z_i\}$ for any ideal triangulation $\bigcup_i\D(z_i)$
of a complete hyperbolic \hbox{3-manifold~$M$} satisfy $\sum_i(z_i)\wedge(1-z_i)=0$.
(This is a 
consequence of the symplectic nature of the NZ relations, as we will explain in more detail 
in the next section.)  Thus to any such triangulation
we can associate a class $\sum_i[z_i]$ in the Bloch group.  But this class is in fact
independent 
of the triangulation, since (modulo some technical points concerning the fact that
the shapes can degenerate to 0 or 1 under 2--3 Pachner moves) any two triangulations are linked
by a series of ``2--3 Pachner moves" in which two tetrahedra sharing a common face are
replaced by the three tetrahedra defined by their two non-shared and two of their three
shared vertices, and the (signed) sum of the shape parameters of these five tetrahedra
is precisely the 5-term relation and does not affect the class of $\sum_i[z_i]$ in the
Bloch group.  Thus one has a class~$[M]\in\B(\C)$. Moreover, from the very definition of
the Bloch group it follows that the function~$D$ extends to a linear map from $\B(\C)$ 
to~$\R$, and from the discussion in the last section we see that the value of $D$ on the
class~$[M]$ is equal to the volume of~$M$. Although we will not use it, we mention that
by a result of Suslin the Bloch group $\B(F)$ of any field~$F$ is isomorphic up to torsion
to the algebraic $K$-group~$K_3(F)$, with $D$ corresponding to the Borel regulator map
from $K_3(\C)$ to $\R$ in the case~$F=\C$.

On the other hand, as described at the end of the last section, the hyperbolic volume
should actually be seen as the imaginary part of a complexified volume taking values
in~$\C/4\pi^2\Z$, so we would like to replace the function~$D(z)$ by some complex-valued
version of the dilogarithm which, even though it may be many-valued at individual
arguments~$z$, becomes one-valued modulo~$4\pi^2$ if we take a linear combination of its
values with arguments belonging to the Bloch group.  This is the idea behind the passage
from the original Bloch group to the extended one.  The first observation
(cf.~\cite{Z:polylogs}) is that the function $L(v):=\Li_2(1-e^v)$ has the derivative 
$v/(e^{-v}-1)$, which is meromorphic and has residues in~$2\pi i\Z$, so that $L$~itself lifts
to a well-defined function from $\C\ssm2\pi i\Z$ to $\C/4\pi^2\Z$ and satisfies the 
functional equation $L(v+2\pi in)=L(v)-2\pi in\log(1-e^v)$ for $n\in\Z$. 
We now introduce the complex 1-manifold
$$
\wh\C \= \bigl\{(u,v)\in\C^2\mid e^u\+e^v\,=\,1\bigr\}\,.
$$
This is an abelian cover of $\C^\times\ssm\{0,1\}$ via $z=e^u=1-e^v$, with Galois
group isomorphic to~$\Z^2$. The extended Bloch group $\hB(\C)$ as defined
in~\cite{Goette} or~\cite{Zickert:extended} is the kernel of the map
$\wh d:\Z[\wh\C]\to\WW(\C)$, where $\WW(\C)$  is defined by requiring only 
$x\wedge y+y\wedge x=0$ (rather than $x\wedge x=0$, which is stronger by 2-torsion)
and where $\wh d$ maps $[u,v]:=[(u,v)]\in\Z[\wh\C]$ to~$u\wedge v$, divided by an appropriate
lifted version of the 5-term relation, namely, the $\Z$-span of the set of elements
$\sum_{j=1}^5(-1)^j[u_j,v_j]$ of $\Z(\wh\C)$ satisfying $(u_2,u_4)=(u_1+u_3,u_3+u_5)$
and $(v_1,v_3,v_5)=(u_5+v_2,v_2+v_4,u_1+v_4)$.  There is an extended
regulator map from $\hB(\C)$ to~$\C/4\pi^2\Z$ given by mapping $\sum[u_j,v_j]$ to
$\sum\calL(u_j,v_j)$, where $\calL(u,v)=L(v)+\frac12uv-\frac{\pi^2}6$, which one
can check vanishes modulo~$4\pi^2$ on the lifted \hbox{5-term} relation.  One can
also define $\hB(F)$ for any subfield~$F$ of~$\C$, such as an embedded number field, by
replacing $\wh\C$ by the subset $\wh F$ consisting of pairs $(u,v)$ with $e^u=1-e^v\in F$.  

As a final remark, one can wonder to what extent studying just hyperbolic 3-manifolds
lets one understand the full Bloch group of~$\Qbar$.  For instance, does every element of
$\B(\Qbar)$ occur as a rational linear combination of the Bloch group invariants of some
hyperbolic 3-manifolds? Even more basically, does every number field with at least one
non-real embedding occur as the trace field of some hyperbolic 3-manifold?  The latter
question was posed explicitly by Walter in~\cite{N:hilbert}.



\section{Symplectic properties}
\label{sec.5}

In retrospect, the symplectic properties as descibed in equation~\eqref{Th2.2} and the
following text, and their refinement from~$\Q$ to~$\Z$ as given in the follow-up
paper~\cite{N:combinatorics},
turned out to be the most important aspects of these papers.  They are responsible both
for the existence of the potential function and for all of the applications
to quantization that we will describe in the next section, as well as many of the connections
to number theory described in Section~\ref{sec.7}.

Define an $N\times2N$ matrix $H=(A\,B)$ whose rows form a $\Z$-basis for the lattice spanned 
by the edge equations~\eqref{edge} together with one ``peripheral'' equation
(a coprime linear combination of the meridian and longitude equations in~\eqref{ML}) at
each cusp. Then the above cited results in~\cite{NZ} imply that $AB^t$ is symmetric and that
$H$ has rank~$N$, meaning that its $2N$ columns generate~$\Q^N$. Together, these two statements
are equivalent to saying that $H$ can be extended to a $2N\times2N$ matrix $\sma ABCD$ in 
$\Sp_{2N}(\Q)$, meaning that ${\sma ABCD}^{-1}=\sma{D^t}{-B^t}{-C^t}{A^t}$.  

But in fact the $2N$ columns of~$H$ span the {\it lattice} $\Z^N$, which is equivalent to 
saying that $H$ can be completed to a $2N\times2N$ symplectic matrix over~$\Z$. 
(We will call such a matrix {\it half-symplectic}.)
This follows from the chain complex defined by Walter in~\cite{N:combinatorics}.  
Explicitly, for any simplex~$\D$, let $J_\D$ as the abelian group generated by $e_1,e_2,e_3$ 
(corresponding to the pairs of opposite edges) subject to the relation $e_1+e_2+e_3=0$. 
This is a free abelian group of rank 2, with a canonical nonsingular, skew-symmetric bilinear 
form given by~\cite[Sec.4]{N:combinatorics}
\be
\label{omega}
\langle e_1, e_2 \rangle \,=\, \langle e_2, e_3 \rangle \,=\, \langle e_3, e_1 \rangle
\,=\, -\langle e_2, e_1 \rangle \,=\, - \langle e_3, e_2 \rangle 
\,=\,- \langle e_1, e_3 \rangle \,=\, 1\,.
\ee
The Neumann chain complex associated to an ideal triangulation is then defined by
\be
\label{ChainComplex}
 0 \;\longrightarrow\; C_0\;\stackrel\a\longrightarrow\; C_1\;\stackrel\b\longrightarrow\; 
 J\;\stackrel{\b^*}\longrightarrow\; C_1\;\stackrel{\a^*}\longrightarrow\;
 C_0\;\longrightarrow 0\,.
\ee
Here $C_0$ and $C_1$ are the free abelian groups on the unoriented 0- and 1-simplices
(cusps and edges), respectively, and $J=\bigoplus_{i=1}^NJ_{\D_i}$ (sum over the 3-simplices
or tetrahedra), while $\a$ maps any cusp the sum of its incident edges, the $J_\D$-component
of $\b$ of any any edge is the sum of the edges of~$\D$ that are identified with it, 
and $\a^*$ and $\b^*$ are the duals of $\a$ and~$\B$ with respect to the obvious scalar
products 
on~$C_i$ and the symplectic form on~$J$.  Walter shows (\cite{N:combinatorics}, Theorem~4.1) 
that the sequence~\eqref{ChainComplex} is a chain complex and, at least after 
tensoring with~$\Z[\frac12]$, is exact except in the middle, where the homology is the
sum of~$h$ rank 2 modules isomorphic to~$H_1(T^2_i)$~$(i=1,\dots,h)$.  Note that the map $\b$
is given precisely by the matrix $R=(R'\,R'')$ as defined in~\eqref{edge} if we choose the
obvious basis for $C_1$ and the basis of $J$ given by choosing the basis $(e_1,e_2)$ for
every~$J_{\D_j}$. The rest of the proof that $H$ is half-symplectic follows easily from the
theorem just quoted and will be left to the reader. 

We make two remarks about this.  The first is that both the construction the chain complex
and the statement about its homology were done in~\cite{N:combinatorics} also for 3-manifolds
with boundary components of arbitrary genus (so the vertices of~$C_0$ need not be cusps),
and of course also do not require any hyperbolic structure.  The other is that the gluing
equations of~\cite{NZ} and the symplectic results of~\cite{N:combinatorics} were extended to
arbitrary $\PGL_n$-representations in~\cite{GZi}.

Half-symplectic matrices occur in other contexts, e.g., in connection with Nahm's
conjecture on the modularity of certain $q$-hypergeometric series, and also lead to a new
description of the Bloch group. Both topics will be discussed in more detail in
Section~\ref{sec.7}.


\section{Quantization}
\label{sec.6}

Perhaps the most far-reaching consequences of Walter's work on the combinatorics of
3-dimensional triangulations have been the applications of the symplectic structure
to quantization.  

Recall the definition of~$J_\D$ for a single tetrahedron~$\D$ as the abelian group
$\langle e_1,e_2,e_3\mid e_1+e_2+e_3=0\rangle$ with the sympectic structure~\eqref{omega}.
This symplectic structure on each space $J_\Delta \otimes_\BZ \BQ$ for any
ideal tetrahedron~$\Delta$ leads to an integral Lagrangian subspace of the
10-dimensional symplectic space \hbox{$\oplus_{j=1}^5 J_{\Delta_j} \otimes_\BZ\BQ$} associated
to 5 tetrahedra $\Delta_1,\dots,\Delta_5$  that participate in a 2--3 Pachner move.
Roughly speaking, the Lagrangian subspace records the linear relations among the
angles of the five tetrahedra, where the signed sum of the angles around each interior 
edge of the Pachner move is zero. 

The quantization of this Lagrangian subspace has appeared numerous times in the
mathematics and physics literature, under different names, and has led to interesting
quantum invariants in dimensions two, three and four. We briefly discuss this now.
In dimension two, Kashaev and independently Fock-Goncharov~\cite{Kashaev:quantization,FG0,FG} 
used the above NZ-symplectic structure to
study the change of coordinates of ideally triangulated surfaces under a 2--2 Pachner
move. They found that the corresponding isomorphism of commutative algebras can be
described in terms of cluster algebras, leading to two dual sets of coordinates
(the so-called $\mathcal{X}$-coordinates and the $\mathcal{A}$-coordinates) whose
quantization leads to a representation of the so-called Ptolemy groupoid, and in 
particular of the mapping class group of a punctured surface, and also of braid groups.
These representations are always infinite-dimensional (because there are no finite square
matrices $A$ and~$B$ satisfying the relation $AB-BA=I$), the Hilbert spaces are typically
$L^2(\BR^n)$ for some~$n$, and the corresponding theory is usually known as quantum 
Teichm\"uller theory.

Going one dimension higher, the $\mathcal{X}$-coordinates of a 3-dimensional ideal
triangulation are nothing but the shapes of the ideal tetrahedra, whereas the
$\mathcal{A}$-coordinates are the Ptolemy variables of the ideal tetrahedra. The latter
are assignments of nonzero complex numbers to the edges of the ideal triangulation
(where identified edges are given the same variable) that satisfy a system of
quadratic equations: a (suitably) signed sum $ab+cd=ef$ where $(a,b)$, $(b,d)$
and $(e,f)$ are the Ptolemy variables of the three pairs of opposite edges. It turns
out that the NZ gluing equations for shapes are equivalent to the Ptolemy equations
(see for instance~\cite{Ga:ptolemy}), and this is not only theoretically interesting, but
practically, too. The quantization of the shape and Ptolemy variables of an ideal
triangulation uses two ingredients, the kinematical kernel of Kashaev~\cite{Kashaev:4D}
and a special
function, the Faddeev quantum dilogarithm that satisfies an integral pentagon identity.
According to Kashaev, the kinematical kernel is nothing but the quantization of the
NZ Lagrangian mentioned above. The outcome of this quantization is the existence of
topological invariants of ideally triangulated 3-manifolds, the invariants being
analytic functions in a cut place $\BC'=\BC\setminus (-\infty,0]$, expressed in terms
of finite dimensional state integrals whose integrand is often determined by
the combinatorial data of an ideal triangulation, namely its Neumann--Zagier matrices.
This construction, that is often known as quantum hyperbolic geometry, has been
axiomatized by Kashaev, and uses as input the combinatorial data of an ideal triangulation
together with a self-dual locally compact abelian group with fixed Gaussian, Fourier
kernel and quantum dilogarithm.   This then leads to further analytic invariants of
3-manifolds, two examples of which are the Kashaev--Luo--Vartanov invariants~\cite{KLV}
and the meromorphic 3D-index~\cite{GK:mero}, for which the LCA groups are $\BR \times \BR$
and $S^1 \times \BZ$, respectively. It is worth noting that the Andersen--Kashaev state
integrals are conjectured to be the partition function of complex Chern--Simons theory
(i.e., Chern--Simons theory with complex gauge group).
The latter is not known to satisfy the cut-and-paste arguments that the $\mathrm{SU}(2)$
Chern--Simons theory does, and as a result, one does not have an a priori definition
of complex Chern--Simons theory other than the state integrals, nor a clear reason
why the infinite dimensional path integral localizes to a finite-dimensional one. 

Finally, going yet one dimension higher, the five ideal tetrehedra that participate in a
2--3 Pachner move form the boundary of a single 4-dimensional simplex, a pentachoron.
(Excuse our Greek.) This gives a 4-dimensional interpretation of the NZ-structure and
of the kinematical kernel, and using a complex root of unity, Kashaev was able to give
a tensor invariant under 4-dimensional Pachner moves and thus construct corresponding
topological invariants of closed, triangulated 4-manifolds at roots of
unity~\cite{Kashaev:4D}. This concludes our discussion of the kinematical kernel in 2, 3
and~4 dimensions.

In a different direction, mathematical physicists, using corresondence principles among
supersymmetric theories, have came up with unexpected constructions of various collections
of \hbox{$q$-series} with integer coefficients associated to 3-manifolds.  Perhaps the most
remarkable of these is the 3D-index of Dimofte--Gaiotto--Gukov \cite{DGG1,DGG2}, where
the $q$-series in question, which in this case are indexed by pairs of integers,
were defined explicitly in terms of the NZ-matrices of a suitable ideal
triangulation, with their coefficients counting the number of BPS states of a supersymmetric
theory. This DGG 3D-index was subsequently shown~\cite{GHRS} to be a topological 
invariant of cusped hyperbolic 3-manifolds, and was also extended to 
a meromorphic function of two variables (in case the boundary of the 3-manifold is a 
single torus) whose Laurent coefficients are the DGG index~\cite{GK:mero}.

A quite different place where the NZ equations appear in quantum topology is in 
connection with the Kashaev invariant and Kashaev's famous Volume Conjecture.  
The Kashaev invariant $\langle K\rangle_n$ is a computable algebraic number that was
defined for any knot~$K$ and any positive integer~$n$ by Kashaev~\cite{K95} in~1995 using
ideas of quantum topology similar to those discussed above, and of which an alternative
definition in terms of the so-called colored Jones polynomial was later found by
H.~and J.~Murakami.  The Volume Conjecture~\cite{K97} says that the logarithm of
$|\langle K\rangle_n|$ is asymptotically equal to $n/2\pi$ times the hyperbolic volume
of the knot complement $M=S^3\ssm K$ whenever $M$ is hyperbolic, a very surprising
connection between hyperbolic geometry and 3-dimensional quantum topology that has
given rise to a great deal of subsequent research and has been refined in many ways
and by several authors in connection with complex Chern-Simons theory. In particular,
one has the conjectural sharpening~\cite{Gar:Viet, DGLZ} 
\be
\label{arithmVolConj}
\langle K\rangle_n \;\, \sim  \;\,
n^{3/2}\, e^{\V(K)n/2\pi i}\,\Phi^K\Bigl(\frac{2\pi i}n\Bigr) 
\ee
to all orders in~$1/n$ as~$n\to\infty$, where $\Phi^K(h)$ is a power series in~$h$ with
algebraic coefficients that can be computed to any order in any explicit example, e.g.
\be
\label{as41}
\Phi^{4_1}(h)\=  \frac1{\sqrt[4]{3}}\, 
\Bigl(1 \+ \frac{11}{72\sqrt{-3}}\,h \+ \frac{697}{2\,(72\sqrt{-3})^2}\,h^2
 \+ \frac{724351}{30\,(72\sqrt{-3})^3}\,h^3 \+\cdots\Bigr)\,. 
\ee
for the $4_1$ (figure~8) knot. In~\cite{DG}, an explicit candidate for this power series 
is constructed for any knot as a formal Gaussian integral whose integrand is defined in
terms of the NZ data of an ideal triangulation of~$M$.  It is not yet known beyond the
leading term that the series constructed there is a topological invariant (i.e.,
independent of the choice of triangulation), although this would of course follow from
the conjecture that the asymptotic formula~\eqref{arithmVolConj} holds with this series.
In a follow-up paper~\cite{DG2}, the construction was extended, still using NZ data in
an essential way, to give explicitly computable power series $\Phi_\a^K(h)$ for
any~$\a\in\Q/\Z$, with $\Phi_0^K=\Phi^K$, that is expected to be the power series
predicted by the quantum modularity conjecture for knots that we will discuss in the
next section.

Finally, it is worth noting that the NZ-equations and their symplectic properties
lead to an explicit quantization of the shape variables, where one replaces each
$z$, $z'$ and $z''$ by operators that suitably commute. This was carried out by 
Dimofte~\cite{Dimofte:quantum}, who defined a quantized version of the gluing equations, 
a so-called quantum curve, which is expected to annihilate the partition function 
of complex Chern--Simons theory and to be ultimately related to the asymptotics of
quantum invariants.


\section{Connections to number theory}
\label{sec.7}

The paper~\cite{NZ} and its sequel~\cite{N:combinatorics} suggested or led to
several interesting developments in pure number theory as well as in topology.
In this final section we describe of a few of these.

\bigskip

\noindent{\bf Values of Dedekind zeta functions and higher Bloch groups}

\smallskip\noindent
An important subclass of hyperbolic manifolds~$M=\BH^3/\G$ are the arithmetic ones, 
where $\G$ is either the Bianchi group $\SL_2(\calO_F)$ for some imaginary quadratic 
field~$F$ or more generally a group of units in a quaternion
algebra over a number field~$F$ of higher degree $n=r_1+2$ having only
one complex embedding up to complex conjugation.  In both cases,
classical results (proved in the first case by Humbert already in~1919) say
that the volume of ~$M$ is a simple multiple (a power of~$\pi$ times the
square-root of the discriminant) of the value at $s=2$ of the Dedekind zeta
function~$\z_F(s)$ of the field~$F$.  An immediate consequence of this and
of the volume formulas discussed in Section~\ref{sec.3} is that this zeta value
is a multiple of a linear combination of values of the Bloch-Wigner dilogarithm
at algebraic arguments.  This consequence was observed in~\cite{Z:hyperbolic}
and was also generalized there to the value of~$\z_F(2)$ for arbitrary number fields~$F$, 
with $[F:\BQ]=r_1+2r_2$ for any value of~$r_2\ge1$.  (If~$r_2=0$, then the well-known 
Klingen-Siegel theorem asserts that $\z_F(2)$ is a rational multiple of
$\pi^{2r_1}\sqrt{D_F}\,$.) Now the group $\SL_2(\calO_F)$ acts as a discrete group of
isometries of $(\BH^2)^{r_1}\times(\BH^3)^{r_2}$ with a quotient of finite volume, and
there are also quaternionic groups~$\G$ over~$F$ that acts freely and discretely
on~$(\BH^3)^{r_2}$, the volume of the quotient in both cases being an elementary
multiple of~$\z_F(2)$. This gives a ``poly-3-hyperbolic'' manifold $M=(\BH^3)^{r_2}/\G$
with volume proportional to~$\z_F(2)$.  A rather amusing lemma says that any such
manifold has a decomposition (disjoint except for the boundaries) into finitely many
$r_2$-fold products of hyperbolic tetrahedra, and it follows that~$\z_F(2)$  for any
number field has an expression as a linear combination of $r_2$-fold products of
values of~$D(z)$ at algebraic arguments, generalizing the Klingen-Siegel theorem
in an unexpected way. 

The connection with 3-dimensional hyperbolic geometry applies only to the values of
Dedekind zeta functions at~$s=2$, but suggested that there might be similar statements for
$\z_F(m)$ with $m>2$ in terms of the $m$th polylogarithm function~$\Li_m(z)$.  Extensive
numerical experiments led to a concrete conjecture saying that this is the case and also to
a definition (originally highly speculative, but now supported by more theory) of ``higher
Bloch groups" $\B_m(F)$ that should be isomorphic after tensoring with~$\Q$ to the higher
algebraic \hbox{$K$-groups} $K_{2m-1}(F)$ and should express the Borel regulator in terms of
polylogarithms. (For a survey, see~\cite{Z:polylogs}.)  This conjecture, now over 30 years
old, has been studied extensively by Beilinson, Deligne, de Jeu, Goncharov, Rudenko and 
others, with the cases~$m=3$ and~$m=4$ now being essentially settled.

\bigskip 

\noindent{\bf Units in cyclotomic extensions of number fields}

\smallskip\noindent
As already mentioned in the last section, the analysis of the Kashaev invariant
and the modular generalization of the Volume Conjecture discussed below led to the
definition of
certain power series $\Phi_\a^K(h)$ associated to a hyperbolic knot and a number $\a\in\Q/\Z$
that can be computed numerically in any given case.  Extensive numerical computations
for simple knots and simple rational numbers~$\a$ suggested that this power series not
only has algebraic coefficients, but that (up to a root of a unity and the square-root
of a number in the trace field~$F_K$ of the knot independent of~$\a$) its $n$th power
belongs to $F_{K,n}[[h]]$, where $n$ is the denominator of~$\a$ and $F_{K,n}=F_K(e^{2\pi i\a})$
the $n$th cyclotomic extension of~$F_K$.  Equivalently, $\Phi_\a^K(h)$ itself is the product
of a power series in $F_{K,n}[[h]]$ with the $n$th root of an element of $F_{K,n}^\times$.
Moreover, in each case studied the latter factor turned out to be the $n$th root of a
{\it unit}, and not just a non-zero number, of~$F_{K,n}$, and in the case of ``sister knots''
(like~$5_2$ and the $(-2,3,7)$-pretzel knot) having the same Bloch group class were the
same for both knots, even though the rest of the power series were completely different.
This led us, together with Frank Calegari, to conjecture and later to prove~\cite{CGZ},
that there was a canonical class of elements in cyclotomic extensions of arbitrary number
fields associated to elements of their Bloch groups, whether or not the fields arise
from topology.  Explicitly, to any number field~$F$ and any element~$\xi$ of the Bloch
group of~$F$ one can associate canonically defined elements of $U(F_n)/U(F_n)^n$ for
every~$n$, where $U(F_n)$ denotes the group of units (more precisely, of $S$-units
for some~$S$ depending on~$\xi$ but independent of~$n$) of the $n$th cyclotomic
extension~$F_n$ of~$F$.  Actually, two quite different constructions were given, one in
terms of an element of~$\B(F)$ and one in terms of an element in~$K_3(F)$, and work in
progress announced in~\cite{CGZ} suggests that there will be a generalization to
$\B_m(F)$ and $K_{2m-1}(F)$ for any~$m>2$.  The construction in terms of the Bloch group
is quite simple, although the proof that it gives units and is independent (up
to~$n$th powers) of all choices is long: if $\xi$ is represented by $\sum[z_i]\in\Z[F]$,
then the number $\prod D_{\z}(z_i^{1/n})$ is the product of an $n$th power in~$F_n$ with
an $S$-unit, and this is the unit we are looking for. Here $\z$ is a primitive $n$th root
of unity and $D_{\z}(x)=\prod_{k=1}^{n-1}(1-\z^kx)^k$ is the ``cyclic quantum dilogarithm''
function, which by a result of Kashaev, Mangazeev and Stroganov satisfies an analogue of
the 5-term relation of the classical dilogarithm.

\bigskip

\noindent{\bf $q$-series and Nahm's conjecture}

\smallskip\noindent
An unexpected consequence of the work on units just described was a proof of one direction
of a conjecture by the mathematical physicist Werner Nahm that had predicted an
extremely surprising connection between the Bloch group and the modularity of certain
$q$-hypergeometric series.  The simplest case of such a ``Nahm sum" is the infinite series
\be
\label{Nahm0}
F_{a,b,c}(q)\= \sum_{n=0}^\infty\frac{q^{\frac12an^2 +bn +c}}{(q)_n}
\qquad(a,b,c\in\Q,\; a>0),
\ee
where $(q)_n=(1-q)(1-q^2)\cdots(1-q^n)$ is the so-called quantum factorial. This function
is known to be modular (in~$\tau$, where $q=e^{2\pi i\tau}$) when $(a,b,c)$ is
$(2,0,-\frac1{60})$ or $(2,1,\frac{11}{60})$ by the famous Rogers-Ramanujan identities
and in a handful of other cases by classical results of Euler, Gauss and others.  It is very
rare that a $q$-hypergeometric series (meaning an infinite sum whose adjacent terms differ by
fixed rational functions of $q$ and $q^n$) is at the same time a modular function, and in fact
for the special series~\eqref{Nahm0} this happens only for seven triples~$(a,b,c)$,
as predicted by Nahm's conjecture and proved in~\cite{Z:dilogarithm}.  Nahm raised the
general question when a (multi-dimensional) $q$-hypergeometric series can be modular and,
motivated by examples coming from characters of vertex operator algebras, discovered a 
possible answer in terms of the Bloch group. Concretely, he generalized~\eqref{Nahm1} to
\be
\label{Nahm1} 
F_{a,b,c}(q)\= \sum_{n_1,\dots,n_N\ge0}
\frac{q^{\frac12n^tan\+b^tn\+c}}{(q)_{n_1}\cdots(q)_{n_N}}
\ee
for any~$N\ge1$, where $a$ is now a positive definite symmetric $N\times N$ matrix with
rational coefficients, $b$ a vector in~$\Q^N$, and $c$ a rational number. There is still no
complete ``if and only if'' conjecture predicting exactly when these Nahm sums are modular
functions, but Nahm gave a precise conjecture for a necessary condition and a partial
conjecture for the sufficiency.  It is the first part that was proved in~\cite{CGZ}, while
the correct formulation and proof of the converse direction is still an active research
subject. 

The modularity criterium that Nahm found depended on his observation that for any
solution $(z_1,\dots,z_N)$ of the system of equations
\be
\label{NahmEq}
 1\m z_i \= \prod_{j=1}^N z_j^{a_{ij}} \qquad(i=1,\dots,N)\,, 
\ee
the element $[z_1]+\cdots+[z_N]$ of~$\Z[\C]$ belongs to~$\B(\C)$. This is a direct consequence
of the symmetry of~$a$, because $\sum_i(z_i)\wedge(1-z_i)=\sum_{i,j}a_{ij}\,(z_i)\wedge(z_j)=0$. 
(This is the same argument as was used in Section~\ref{sec.5} for the corresponding statement
for the Neumann-Zagier equations, and applies more generally to all half-symplectic matrices,
as discussed below.)  The ``only if'' direction of Nahm's conjecture then says that
$F_{a,b,c}(q)$ can be modular only if this element of the Bloch group vanishes for the
unique solution of the Nahm equation having all $z_i\in(0,1)$.  The proof, given
in~\cite{CGZ}, uses both the results there about the units coming from Bloch elements
as described above and an asymptotic analysis of Nahm sums near roots of unity published
separately by the two of us.

\bigskip

\noindent{\bf Half-symplectic matrices and the Bloch group}

\smallskip\noindent
At the end of Section~\ref{sec.5}, we saw how the NZ equations lead to a ``half-symplectic
matrix,'' meaning the upper half $H=(A\,B)$ of a $2N\times2N$ symplectic matrix $\vsma ABCD$ 
over~$\Z$. To any such matrix we associate the system of polynomial equations
\be
\label{NZeq}
\prod_{j=1}^N z_j^{A_{ij}} \=
(-1)^{(AB^t)_{ii}\vphantom{y_{y_y}}}\,\prod_{j=1}^N (1-z_j)^{B_{ij}}
\qquad(i=1,\dots,N)\,.
\ee
This is a straight generalization of the original NZ equations in the topological setting
except possibly for the sign, but we checked in thousands of examples using Snappy~\cite{snappy}
that this is the right sign, and this could presumably be proved using the ``parity 
condition'' in~\cite{N:combinatorics}.  Another special case of~\eqref{NZeq} is the Nahm 
equation~\eqref{NahmEq}, at least when the matrix~$a$ is integral and even, the
half-symplectic matrix then being~$(1\,a)$ and the full symplectic matrix~$\vsma 1a01$.

For any solution~$z$ of~\eqref{NZeq}, the element $2\sum_{j=1}^N[z_j]$ belongs to the usual
Bloch group by the same argument as was used for Nahm sums. We can in fact divide by~2
and lift to an element $[H,z]$ of the extended Bloch group by setting
\be
\label{Hz}
[H,z] \= \sum_{j=1}^N [u_j,v_j] \+ [\xi,\xi'] \+ [-\xi,\,\xi'-\xi+\pi i]\,,
\ee
where $u_j$ and $v_j$ are arbitrary choices of logarithms of $z_j$ and~$1-z_j$, respectively, 
$\xi$ is defined as $\frac1\pi\,(Au-Bv)^t(Cu-Dv)$ for any completion $(C\,D)$ of~$H$ to a full
symplectic matrix, and $\xi'$ is any choice of logarithm of~$1-e^\xi$.
The facts that this is in the kernel of~$d:\Z[\wh\C]\to\WW(\C)$ and that its image modulo
extended 5-term relations is independent of all choices (at least modulo 8-torsion in the
Bloch group of the field generated by the~$z$'s) can be checked by direct computations which
are sketched in Section~6.1 of~\cite{GZ:kashaev}, together with a more precise form that
eliminates the torsion ambiguity.

We have the five following equivalence relations among pairs, each motivated by a change of
the choices made in the topological situation, that do not change this class:
\begin{itemize}
\item Stability: increase~$N$ by $N+1$, replace $A$ and~$B$ in~$H=(A\,B)$ by their direct
sums with $(1)$ and~$(0)$, respectively, and set $z_{N+1}=1$, corresponding in the topological
case to adding a degenerate simplex to a 3-manifold triangulation.
\item Changing the equations: multiply $H$ on the left by an element of~$\GL_N(\Z)$ without
  changing the $z$'s.  This corresponds to replacing the $N$ relations~\eqref{NZeq} by
  multiplicative combinations of them in an invertible way.
\item  Renumbering: multiply~$H$ on the right  by an $N\times N$ permutation matrix, and
  permute the $z_j$'s by the same matrix, corresponding to a renumbering of the simplices of
  a triangulation. 
\item  New shape parameters: for each $j=1,\dots,N$, multiply the $N\times2$ matrix
  $(A_j\,B_j)$ (where $A_j$ and $B_j$ denote the $j$th columns of $A$ and~$B$, respectively)
  by a power of the element $\sma0{-1}1{-1}$ of order~3 in~$\SL_2(\Z)$, and replace $z_j$
  correspondingly by $z_j$, $z_j'=1/(1-z_j)$, or $z_j''=1/(1-z_j)$.
\item Algebraic 2--3 Pachner moves: remove two columns of $A$ and the same two columns of
  $B$ and replace them by 3 new columns that are specific $\Z$-linear combinations,
  simultaneously replacing $N$ by $N+1$ and changing two of the $z_j$ to three others in
  such a way that the corresponding change of $[H,z]$ is a 5-term relation.
  The explicit formulas were first written down in the special case corresponding to the
  Nahm sums~\eqref{Nahm1} by Sander Zwegers in an unpublished 2011 conference talk and were
  then given for abritrary symplectic matrices in Equation (3-27) of~\cite{DG}. This
  corresponds to stabilizing $H$ three times and then multiplying it on the left by a specific 
  element of~$\Sp_{10}(\Z)$, and then unstabilizing three times.
\end{itemize}
This gives us a new abelian group that maps to the extended Bloch group, namely the set of
all pairs $(H,z)$ as above modulo these equivalence relations, with addition given by
direct sum. In fact this map is an isomorphism, meaning that any element of the Bloch group
can be realized by some half-symplectic matrix and solution of the coresponding generalized
Neumann-Zagier equations and that any 5-term relation can be lifted to one coming from an
algebraic 2--3 Pachner move. A more complete discussion is given in Section~6.1
of~\cite{GZ:kashaev} (full details will be given later), while Section~6.3 of the same paper
shows how to attach Nahm sum-like $q$-series to arbitrary half-symplectic matrices.

\bigskip 

\noindent{\bf From the Kashaev invariant to quantum modular forms}

\smallskip\noindent 
Nahm's conjecture already highlighted a connection between half-symplectic matrices and 
questions of modularity, but there are other and more direct connections between hyperbolic 
3-manifolds and the modular group~$\SL_2(\Z)$ that we now describe.

At the end of Section~\ref{sec.6} we discussed Kashaev's volume conjecture and its 
refinement~\eqref{arithmVolConj}.  That statement in turn was generalized in~\cite{Z:QMF}
on the basis of numerical computations to a conjectural asymptotic formula having a strong
modular flavor. To state it, we first note that the Kashaev invariant $\la K\ra_n$ of a
knot~$K$ can be generalized to a function $J^K:\Q/\Z\to\Qbar$ whose value at $-1/n$ for any
$n\ge1$ is $\la K\ra_n$ and which is $\text{Gal}(\Qbar/\Q)$-equivariant.
Then the conjectural generalization of~\eqref{arithmVolConj} is the statement that
\be
\label{QMCb}
J^K\Bigl(\frac{an+b}{cn+d}\Bigr) \;\sim\; (cn+d)^{3/2} \;
e^{\Vol_\C(K)(n+d/c)/2\pi i}\; \Phi^K_{a/c}\Bigl(\frac{2 \pi i}{c(cn+d)}\Bigr)\; J^K(n)
\ee
for every matrix $\sma abcd\in\SL_2(\Z)$ as $n$ tends to infinity through either 
integers or rational numbers with bounded denominator, with~\eqref{arithmVolConj} being
the special case when $\sma abcd=\sma 0{-1}10$ and $n$ is integral.  Here $\Phi^K_\a(h)$
for $\a\in\Q/\Z$ is a power series that is conjectured to be the one constructed in~\cite{DG2}
and discussed at the end of Section~\ref{sec.6}.

In~\cite{GZ:kashaev}, this modularity conjecture was verified experimentally for a few knots 
to many terms and to a high degree of precision, and was also successively refined in
several different directions, the final statement being the existence of a whole matrix
$\J^K$ of $\Qbar$-valued functions on~$\Q/\Z$ (generalized Kashaev invariants) which
conjecturally has much better modularity properties than the original scalar function~$J^K$.
Explicitly, \eqref{QMCb} lifts to a similar statement with $J^K$ replaced by the
matrix~$\J^K$ and the completed formal power series
$\wh\Phi^K_{a/c}(h)=e^{\Vol_\C(K)/2\pi ic^2h}\,\Phi^K_{a/c}(h)$ by a matrix $\wh\bPhi^K_{a/c}(h)$
of completed formal power series acting by right multiplication. 
The rows and columns of these matrices are indexed by the bounded parabolic
flat connections, or equivalently by an index~0 (trivial connection) and indices~$1,\dots,r$
corresponding to the solutions of the NZ equation for a triangulation of the knot complement,
with the original scalar-valued functions $J^K$ and $\Phi^K$ being the $(0,1)$ and $(1,1)$ 
entries of $\J^K$ and $\bPhi^K$, respectively.  The really new aspect is that, by replacing
the original scalar functions by matrices, we obtain a matrix of completed formal power
series in~$h$ that (conjecturally, like everything else in this story) extend to
real-analytic functions on the 
positive and negative real line and in fact to holomorphic functions on the two cut planes
$\C\ssm(-\infty,0]$ and $\C\ssm[0,\infty)$. This discovery, which arises through the
possibility of associating matrix-valued $q$-series to the knot~\cite{GZ:qseries}, gives
rise to the new concept of ``holomorphic quantum modular form" that then turned out to
appear also in many other situations, including various known modular objects like mock
modular forms or Eisenstein series of odd weight on the full modular group.

Finally, we mention that the new generalized Kashaev invariants have beautiful arithmetic
properties generalizing the known property~\cite{Hab} that the original Kashaev invariant
belongs to the so-called Habiro ring~$\calH=\varprojlim\Z[q]/(q)_n\,$.  The results
of~\cite{GZ:kashaev} and~\cite{GZ:qseries} suggest that there should be a Habiro
ring~$\calH_F$ associated to any number field~$F$ in which the generalized Kashaev invariants
take their values and which is graded by the Bloch group of~$F$. (The latter property is
invisible in the classical case since $\B(\Q)\otimes\Q=\{0\}$.)  We are currently working
on this jointly with Peter Scholze, and already have a candidate for~$\calH_F$, as well as
a partial lifting of the algebraic units of~\cite{CGZ} to formal power series with
Habiro-like properties.

\medskip 


\bibliographystyle{plain}
\bibliography{biblio}
\end{document}